\documentclass[12pt]{amsart}
\usepackage[english]{babel}
\usepackage[utf8]{inputenc}

\usepackage{fullpage,multicol}

\usepackage{tabu}

\usepackage{amsmath}
\usepackage{amssymb,mathrsfs,enumerate}

\numberwithin{equation}{section}

\newtheorem{thm}{Theorem}[section]
\newtheorem{defi}[thm]{Definition}

\newtheorem{con}[thm]{Conjecture}
\newtheorem{ques}[thm]{Question}

\def\bC{\mathbb C}
\def\pa{\partial}
\def\sO{\mathscr O}
\def\sD{\mathscr D}
\def\al{\alpha}
\def\bZ{\mathbb Z}
\def\bN{\mathbb N}
\def\be{\begin{equation}}
\def\ee{\end{equation}}
\def\Exp{\textup{Exp}}
\def\ba{\mathbf m}

\def\bbf{\mathbf {f}}
\def\bbe{\mathbf e}
\def\sP{\mathscr P}
\def\sC{\mathscr C}

\begin{document}


\title{Bernstein-Sato ideals}

\author{Nero Budur}
\address{Department of Mathematics, KU Leuven, Celestijnenlaan 200B, 3001 Leuven, Belgium;  Yau Mathematical Sciences Center, Tsinghua University, 100084 Beijing, China;  Basque Center for Applied Mathematics, Mazarredo 14, 48009 Bilbao, Spain.}
\email{nero.budur@kuleuven.be}

\author{Robin van der Veer}
\email{robinvdveer@gmail.com}

\author{Lei Wu} 
\address{School of Mathematical Sciences, Zhejiang University,
 Zijingang Campus, Hainayuan 915, Hangzhou, China.}
\email{leiwu23@zju.edu.cn}

\author{Peng Zhou}
\address{Department of Mathematics, University of California,
 753 Evans Hall, Berkeley CA, 94720, USA.}
\email{pzhou.math@gmail.com}

\begin{abstract}
In this paper, we review several results on the zero loci of Bernstein-Sato ideals related to singularities of hypersurfaces.
\end{abstract}

\maketitle

\setcounter{tocdepth}{1}
\tableofcontents

\section{Introduction}

This is an exposition for the Frontiers of Science Awards in Mathematics presenting  results from our article \cite{BS1}, with history, motivation, and further developments.

\subsection{Singularities}  Around 1970, without being precise about the timing, there were two separate developments relevant for us: one in singularity theory, the other in special functions theory.  The first is about the topology of a hypersurface singularity defined by the zero locus of the germ $$f:(\bC^n,0)\to (\bC,0)$$ of a holomorphic function. Milnor, for isolated singularities, and then later Hamm-Le for non-isolated singularities of $f^{-1}(0)$, showed that the restriction of $f$ to the complement of $f^{-1}(0)$ in a small ball centered at $0\in\bC^n$ is a smooth locally trivial fibration onto a small punctured disc. The deviation from global triviality of the fibration is captured by the monodromy action on the cohomology of the Milnor fiber, induced by going once counter-clockwise in a small loop in the disc. The monodromy theorem of Brieskorn, Landman, Grothendieck, then states that the eigenvalues of the monodromy are always roots of unity. Later with the introduction of derived categories and functors, these topological invariants were encoded in one single object by Deligne, the nearby cycles complex  $\psi_f\bC_X$, for $X$ a small analytic open neighborhood of $0\in\bC^n$, and $\bC_X$ the constant sheaf on $X$. The nearby cycles complex comes equipped with a monodromy action, and the cohomology of its stalks at $x\in f^{-1}(0)$ is that of the Milnor fiber of $f$ at $x$ equipped with the monodromy. The union of all sets of eigenvalues of the monodromy for $x\in X$ varying along the zero locus of $f$ is called the {\it monodromy support} $S_f\subset\bC^*$, and is thus a set consisting of roots of unity.

\subsection{$\sD$-modules}

The second development came from two directions, one from algebraic analysis initiated by Sato to use sheaf theory to study hyperfunctions and microfunctions, and one from Bernstein initiated to answer a question of Gelfand about the meromorphic extension in $s$ of complex powers $f^s$. Sato, in certain special cases, and Bernstein in the case $f$ is a polynomial, discovered that there is always a non-zero polynomial $b(s)\in\bC[s]$ such that 
$$
b(s)f^s=P\cdot f^{s+1}
$$
for some linear partial differential operator $P$ on $X=\bC^n$ with polynomial coefficients in $s$, that is, $P\in D_X[s]$, where $D_X=\bC[x_1,\ldots,x_n]\langle \pa_1,\ldots,\pa_n\rangle$, and $\pa_i=\pa/\pa x_i$. Here the differentiation in the right-hand side is formally defined as usual, namely $\pa_i\cdot f^{s}=sf^{-1}\pa_i(f)f^s$ in $\bC[x_1,\ldots,x_n][f^{-1}]\cdot f^s$.
The monic generator of the ideal in $\bC[s]$ generated by such $b(s)$ is called  the {\it Bernstein-Sato polynomial}, or the {\it $b$-function} of $f$, denoted $b_f(s)$.  Bj\"ork extended  this result to the analytic germ case, which we refer to as the local case. The lowest common multiple of the local $b$-functions of the germs of a non-constant regular function $f:X\to \bC$ at $x$ varying in $f^{-1}(0)$, is the $b$-function $b_f(s)$, for any smooth complex affine variety $X$. 

This result marked the birth of $D$-module theory, where it is a fundamental step that $\sO_X[f^{-1}]$ is a finitely generated  $\sD_X$-module. Bernstein used his theorem and integration by parts to show that $f^s$ admits a meromorphic continuation to the whole complex plane with the set of poles equal to the shifts by $\bZ_{\le 0}$ of the set of roots of the $b$-function $b_f(s)$, answering Gelfand's question.  Kashiwara showed that the roots of $b_f(s)$ are negative rational numbers. Malgrange, and then Kashiwara in a more general case, showed  that the two developments above are very much related, namely that by multiplying the roots of the $b$-function of $f$ by $2\pi i$ and then exponentiating, one obtains the monodromy support of $f$. Much later, in 1990s, Oaku used Gr\"obner bases for non-commutative rings such as $D_X$ to give the first algorithm for computing $b_f(s)$, hence for computing the monodromy support $S_f$ as well.

\subsection{}
We summarize the above results in:

\begin{thm}\label{thm1}  Let $f:X\to \bC$ be a regular function on a smooth complex affine irreducible algebraic variety, or the germ at $x\in X$ of a holomorphic function on a complex manifold, such that $f$ is not invertible. Then:

\begin{enumerate}[(i)]
\item (Malgrange, Kashiwara)  There is an equality
$$
\{exp(2\pi i\al)\mid \al\text{ is a root of }b_f(s)\}  = S_f
$$
between the set of exponentials of the roots of the Bernstein-Sato polynomial of $f$ and the monodromy support of $f$.

\item (Kashiwara) The roots of $b_f(s)$ are negative rational numbers.

\item (Monodromy Theorem) The monodromy eigenvalues on the nearby cycles complex of $f$ are roots of unity.

\end{enumerate}
\end{thm}

For precise references we refer to the introduction of \cite{BS1}.

As for other explicit instances of the Riemann-Hilbert correspondence, the topological invariant $S_f$ of $f$ is thus computable via the algebraic/analytic invariant $b_f$. It is well-known that the Bernstein-Sato polynomial $b_f$ is not a topological invariant of $f$, see \cite{Y} where a compendium of examples was calculated by hand, and, in particular, examples are given showing that $b_f$ can change in  deformations of isolated hypersurface singularities within the same embedded-topological equivalence class.

On the other hand, some roots of $b_f$ may well be topological. The negative of the maximal root of $b_f$ is called the {\it log canonical threshold} of $f$.

\begin{ques} (\cite{B-s}) Is the log canonical threshold an embedded-topological invariant for germs of reduced hypersurface singularities?
\end{ques}

The log canonical threshold does not change for deformations of isolated hypersurface singularities within the same embedded-topological equivalence class, by an older result of Varchenko. Beyond this, the closest is a result of \cite{Mc} who showed using Floer theory that for isolated hypersurface singularities the log canonical threshold is an embedded-contact invariant of the link. For topological roots of $b_f$ in the case of plane curves, see \cite{Bl}.

The article \cite{BS1} concerns the generalization of  Theorem \ref{thm1} to the case when $f$ is replaced by a tuple $f_1,\ldots,f_r$.

\section{Motivation} 

The motivation to study a morphism
$$
F=(f_1,\ldots,f_r) : X\to \bC^r
$$
came in fact from a different angle: the monodromy  conjecture of Igusa, Denef-Loeser. 

\subsection{Monodromy conjecture}

Atiyah had answered Gelfand's question on complex powers $f^s$ of a polynomial $f\in\bC[x_1,\ldots,x_n]$ using a different method than Bernstein, namely Hironaka's resolution of singularities. A similar question in the $p$-adic setup, for the so-called $p$-adic local zeta function, was shown by Igusa to be equivalent to the conjecture posed by Borevich-Shafarevich that the generating series over $m\in\bN$ for the number of solutions of $f$ modulo prime powers $p^m$ (say $f$ is defined over $\bZ$ here) is a rational function. Using resolution of singularities, Igusa then proved this conjecture. Based on the analogy with Gelfand's question and a few worked out examples, Igusa then asked if it is always the case that the real parts of the poles of the $p$-adic local zeta function of $f$ are roots of the $b$-function of $f$, for almost all primes $p$.

Denef and Loeser defined the motivic zeta function $Z^{mot}_f(s)$ of a polynomial $f$ in $\bC[x_1,\ldots,x_n]$ as, essentially, a generating series over $m\in\bN$ for the classes in the Grothendieck ring of complex varieties of contact loci $\mathcal X_m(f)$, where $\mathcal X_m(f)$ is the set of algebraic $m$-jets on $\bC^n$ with contact order along $f$ equal to $m$. They showed that the motivic zeta function of $f$ is rational in a certain sense and that it specializes to the $p$-adic local zeta function of Igusa for almost all primes $p$. They conjectured further that the poles of the motivic zeta function are roots of the $b$-function of $f$, which would imply a positive answer to Igusa's question for the $p$-adic local zeta functions.

The monodromy conjecture is still open. The weaker version, implied by the result of Malgrange and Kashiwara, that the exponentials of the poles of the motivic zeta function should give monodromy eigenvalues of $f$, is  open even for isolated singularities.

 \subsection{Generalization}

The following exercise in positive thinking, along with some computer calculations, was done in \cite{B} and it led to some conjectures: suppose that the monodromy conjecture is true, then some version of it should also be true if $f$ is replaced by a tuple $F=(f_1,\ldots,f_r)$ of polynomials $f_j\in\bC[x_1,\ldots,x_n]$. In fact, two versions of it should be true: one relating a multi-variable motivic zeta function $Z^{mot}_F(s_1,\ldots,s_r)$ with a monodromy support locus for $F$ in $(\bC^*)^r$, and, secondly, with a Bernstein-Sato type ideal in $\bC[s_1,\ldots,s_r]$. Moreover, a generalized form of the result of Malgrange and Kashiwara should imply that the version in terms of Bernstein-Sato ideals is stronger than the version in terms of the monodromy support of $F$. Also, one hopes that this generalization provides more elbow room for ideas to prove the monodromy conjecture.

The generalized version of the monodromy conjecture in terms of the monodromy support had already been posed by Loeser. For the motivic zeta function $Z^{mot}_F(s_1,\ldots,s_r)$ one uses as definition the easy generalization of a formula of Denef-Loeser for $Z_f^{mot}(s)$ in terms of an embedded resolution of $f=\prod_{j=1}^rf_j$. The independence of the choice of resolution is proved using the weak factorization theorem. From the definition one has that $Z_F^{mot}(s_1,\ldots,s_r)$ is a rational function on $\bC^r$ with denominator a product of polynomials of type  
\be\label{eq1}
a_1s_1+\ldots+a_rs_r+b,\quad \text{with}\quad a_1,\ldots,a_r,b\in\bN, \text{ and } b\neq 0.
\ee
The union of these hyperplanes in $\bC^r$ is called the {\it polar locus} of the motivic zeta function of $F$.

For the definition of the monodromy support $S_F$ in $(\bC^*)^r$  one uses as definition an easy generalization of the nearby cycles complex of $f$, denoted similarly $\psi_F\bC_X$, and used for the first time by Sabbah. Instead of only one monodromy action, $\psi_F\bC_X$ has $r$ simultaneous monodromy actions, namely going once around each hypersurface $f_j=0$ for $j=1,\ldots, r$. Hence the support of $\psi_F\bC_X$, viewed as a complex of $\bC[t_1^\pm,\ldots,t_r^\pm]$-modules whose cohomology  stalks  are finitely generated, is a closed subscheme of $(\bC^*)^r$. 

\begin{defi}
We call the {\it monodromy support $S_F$} in $(\bC^*)^r$ the underlying reduced closed subscheme of the $\bC[t_1^\pm,\ldots,t_r^\pm]$-support of the complex $\psi_F\bC_X$. 
\end{defi}
Unlike for the motivic zeta function,  the monodromy support can be defined in the analytic setup as well, generalizing that for a single germ of a holomorphic function $f$.

On the other hand there are infinitely many types of Bernstein-Sato ideals for $F$ which generalize the Bernstein-Sato polynomial $b_f(s)$ of a single function $f$. Calculations singled out the following ideal.

\begin{defi}
Define {\it the Bernstein-Sato ideal of F} the ideal $B_F$ generated by all polynomials $b\in\bC[s_1,\ldots,s_r]$ such that there exists $P\in D_X[s_1,\ldots,s_r]$ with
$$
b(s_1,\ldots,s_r)f_1^{s_1}\ldots f_r^{s_r}=P\cdot f_1^{s_1+1}\ldots f_r^{s_r+1}.
$$
\end{defi}
As before, this can be done in the analytic setup as well.

 The generalized versions of the monodromy conjecture are as follows, defining coordinate-wise the map $$\Exp:\bC^r\to(\bC^*)^r,\quad \al\mapsto \exp(2\pi i\al).$$

\begin{con}\label{conj1} Let $F=(f_1,\ldots,f_r)$ with $f_j\in\bC[x_1,\ldots,x_n]$ such that $f=\prod_{j=1}^rf_j$ in not constant.
\begin{enumerate}[(i)]
\item (Loeser) The image under $\Exp$ of the polar locus of the motivic zeta function of $F$ is the monodromy support $S_F$ of $F$.
\item (\cite{B}) The polar locus of the motivic zeta function of $F$  is included in the zero locus of the Bernstein-Sato ideal $B_F$. 
\end{enumerate}
\end{con}

Futhermore the following were posed in \cite{B}:

\begin{con}\label{conj2} 
Let $F=(f_1,\ldots,f_r):X\to\bC^r$ be a morphism of smooth complex affine irreducible algebraic varieties, or the germ at $x\in X$ of a holomorphic map on a complex manifold, such that not all $f_i$ are invertible. 
\begin{enumerate}[(i)]
\item The Bernstein-Sato ideal $B_F$ is generated by products of linear polynomials of type (\ref{eq1}).
\item The image under $\Exp$ of the zero locus of the Bernstein-Sato ideal $B_F$ is the monodromy support locus $S_F$.
\end{enumerate}
\end{con}

 Part $(i)$ of this conjecture is still open.
 
 \section{Results}
 
 The main result of \cite{BS1} is the proof of the second part of Conjecture \ref{conj2} generalizing the classical result of Malgrange and Kashiwara. This completed the following generalization of Theorem \ref{thm1} from the case of a single function $f$ to the case of a tuple of functions $f_1,\ldots, f_r$:

\begin{thm}\label{thrmE} Let $F=(f_1,\ldots,f_r):X\to\bC^r$ be a morphism of smooth complex affine irreducible algebraic varieties, or the germ at $x\in X$ of a holomorphic map on a complex manifold, such that not all $f_i$ are invertible. Let $Z(B_F)$ be the zero locus in $\bC^r$ of the Bernstein-Sato ideal of $F$. Then:

\begin{enumerate}[(i)]
\item $\textup{Exp}(Z(B_F))=S_F.$ 

\item 
\begin{enumerate}[(a)]
\item Every irreducible component of $Z(B_F)$ of codimension 1 is a  hyperplane of type (\ref{eq1}). 

\item Every irreducible component of  $Z(B_F)$ of codimension $>1$  can be translated by an element of $\bZ^r$ inside a component of codimension 1.
\end{enumerate}
\item $S_F$ is a finite union of torsion-translated complex affine subtori of codimension 1 in $(\bC^*)^r$.
\end{enumerate}
\end{thm}

Thus in the algebraic case, the Bernstein-Sato ideal $B_F$  provides an algebraic computation, readily implemented in computer algebra systems, of the monodromy support $S_F$, an embedded topological invariant of the tuple $F=(f_1,\ldots,f_r)$.

It took some years for the statement of this theorem to be reached. The first assertion of part $(ii)$, about the components of codimension 1 of $Z(B_F)$, is much older and proved by  Sabbah in 1987, and Gyoja in 1990 for the condition $b\neq 0$ in (\ref{eq1}),  in the course of showing the existence of a non-zero element in the Bernstein-Sato ideal. As mentioned before, for the missing precise references in this note see the introduction of \cite{BS1}.

Chronologically next was a partial result regarding part $(iii)$ by Sabbah in 1990 who showed that the monodromy support $S_F$ is included in a finite union of torsion-translated complex affine subtori of codimension 1 of $(\bC^*)^r$. Here a complex affine subtorus  means an algebraic subgroup $G\subset (\bC^*)^r$ such that $G\cong (\bC^*)^p$ as algebraic groups for some $0\le p\le r$.

Then the first author and B. Wang proved in 2017 that every irreducible component of $S_F$ is a torsion-translated complex affine subtorus. Nowadays this is part of a more general result on absolute subsets of the space of rank-1 local systems on smooth complex algebraic varieties generalizing the monodromy theorem and various previous results on cohomology jump loci. For this, the monodromy support locus $S_F$ was reinterpreted as the locus of rank-1 local systems $L$ on $(\bC^*)^r$ for which the pullback local system $F_U^{-1}L$ has non-trivial cohomology in the complement  in some small ball in $X$ centered at some point along $f^{-1}(0)$, where $F_U:U=X\setminus f^{-1}(0)\to(\bC^*)^r$ is the restriction of $F:X\to \bC^r$,  and $f=f_1\ldots f_r$.

The last missing part of $(iii)$, namely that all the irreducible components of the monodromy support locus are of codimension 1, was proven in \cite{BLSW}. The proof uses firstly the property of the nearby cycles complex $\psi_f\bC_X$, stemming from the perversity of a shift of this complex, that every eigenvalue in its monodromy support is a zero or a pole of the monodromy zeta function of $f$ at some point along $f^{-1}(0)$. Secondly, the proof uses the compatibility of  the monodromy support with respect to specializations of $F$, that is, the restriction of $S_F$ to 1-parameter subgroups $\bC^*\to (\bC^*)^r$, $\lambda\to (\lambda^{m_1},\ldots,\lambda^{m_r})$ with $m_j\in\bZ_{>0}$ is the monodromy support $S_g$ of $g=f_1^{m_1}\ldots f_r^{m_r}$. This compatibility had been shown earlier by Sabbah.

In light of part $(i)$, which was still a conjecture at the time, it was thus a subconjecture that part $(iii)$ would hold for $\Exp(Z(B_F))$ in place of $S_F$. This is  equivalent to the second assertion in part $(ii)$, about the smaller-dimensional components of $Z(B_F)$ disappearing inside the codimension 1 components after translation by some element of $\bZ^r$. This subconjecture was confirmed unconditionally by Maisonobe \cite{M} using homological algebra over the non-commutative ring $\sD_X[s_1,\ldots,s_r]$.

Part $(i)$  was conjectured in \cite{B}, where one inclusion was also proved, namely that $\textup{Exp}(Z(B_F))$ contains $S_F$. For the proof  in  \cite{BS1} of part $(i)$, all the previous results including the result of Maisonobe  play a crucial role. A major hurdle is that the compatibility of forming the zero locus of the Bernstein-Sato ideal with specializations of $F$ to $g$ as above is still an open question, and in fact, one might expect that there are counterexamples. However, for  specializations generic in a certain sense this compatibility holds. The interplay of this and homological algebra for relative holonomic $\sD$-modules reduces the proof of part $(i)$ to the case when $r=1$, namely, to the classical result of Malgrange and Kashiwara.

Not long after \cite{BS1} a different proof of part $(i)$ was given by R. van der Veer \cite{V} using parametric Gr\"obner bases for relative holonomic $\sD$-modules. His proof does not rely on \cite{M} and in fact it reproves the result of Maisonobe using \cite{BLSW}. Thus \cite{V} also gives a different proof of the second part of $(ii)$.

We note that the following subconjecture of Conjecture \ref{conj2} $(i)$ is still open: 

\begin{ques}
Are the irreducible components of codimension $>1$ of the zero locus of the Bernstein-Sato ideal $B_F$ also linear subvarieties of $\bC^r$?
\end{ques}

\section{Generalization}

There are infinitely many ideals of Bernstein-Sato type one can attach to a tuple $F=(f_1,\ldots,f_r)$ of polynomials or germs of holomorphic functions. 

\begin{defi}
Taking  $$\ba=(m_1,\ldots, m_r)\in\bN^r$$ let $$
B_F^{\,\ba}\subset\bC[s_1,\ldots,s_r]$$ be the ideal generated by all polynomials $b(s_1,\ldots,s_r)$ such that $$b(s_1,\ldots,s_r)\prod_{i=1}^rf_i^{s_i}=P\cdot\prod_{i=1}^rf_i^{s_i+m_i}\text{ for some }P\in D_X[s_1,\ldots,s_r].
$$
Denote its zero locus by $$Z(B_F^{\,\ba})\subset\bC^r.$$ 
\end{defi}
The structure and the precise topological information contained by $Z(B_F^{\,\ba})$ was analysed in \cite{BS2}.

Note that if $\ba=(1,\ldots,1)$ then one obtains the Bernstein-Sato ideal $B_F$ from the previous theorem. Although there are infinitely many ideals $B_F^{\,\ba}$, after taking the image under $\Exp$ of their zero loci there are only finitely many sets to consider, since by \cite{B},
$$
\Exp(Z(B_F^{\,\ba})) = \bigcup_{i}\Exp (Z(B_{F}^{\,\bbe_i})) 
$$
where the union is taken over all $1\le i\le r$ with $m_i\ne 0$ and $f_i$ not invertible, and $\bbe_i$ is the  standard basis vector $(0,\ldots,0,1,0,\ldots, 0)\in\bN^r$ with the unit on the $i$-th position.

On the topological side one has the following:

\begin{defi} Assuming that $\bbf^\ba=\prod_{i=1}^rf_i^{m_i}$ is not invertible,  define 
$$
S_F^\ba  \subset (\bC^*)^r
$$
to be the support, taken with the reduced structure, of the monodromy action on the restriction of $\psi_F\bC_X$ to the zero locus of $\bbf^\ba$. 
\end{defi}
Hence by definition
$$
S_F^\ba= \bigcup_{i}S_F^{\bbe_i}
$$ 
where the union is taken over all $1\le i\le r$ with $m_i\ne 0$ and $f_i$ not invertible, and $S_F^{\bbe_i}$ is the monodromy support of the restriction of $\psi_F\bC_X$ to the zero locus of $f_i$.

\begin{thm}\label{thm3}
Let $F=(f_1,\ldots,f_r):X\to\bC^r$ be a morphism of smooth complex affine irreducible algebraic varieties, or the germ at $x\in X$ of a holomorphic map on a complex manifold.  Let $\ba\in\bN^r$ such that $\bbf^\ba=\prod_{i=1}^rf_i^{m_i}$ is not invertible. Then: 
\begin{enumerate}[(i)]
\item $\textup{Exp} (Z(B_F^{\,\ba})) = S_F^\ba.$
\item 
	\begin{enumerate}[(a)]
		\item Every irreducible component of $Z(B_F^{\,\ba})$ of codimension 1 is a  hyperplane of type 			(\ref{eq1}) and for each such hyperplane there exists $i$ with $m_i\ne 0$ such that $a_i>0$ in (\ref{eq1}). 
		\item Every irreducible component of  $Z(B_F^{\,\ba})$ of codimension $>1$  can be translated by an 		element of $\bZ^r$ inside a component of codimension 1.
	\end{enumerate}
\item $S_F^\ba$ is a finite union of torsion-translated complex affine algebraic subtori of codimension 1 in $(\bC^*)^r$.
\end{enumerate}
\end{thm}

Part $(ii) (a)$ is due to Sabbah and Gyoja. The rest was proven in \cite{BS2} by combining and extending   arguments from \cite{B, BLSW, M, BS1, WZ}.

\section{Other developments}

\subsection{Riemann-Hilbert correspondence for $\psi_F\sP$.}
An explicit Riemann-Hilbert correspondence for Sabbah's  generalization to a tuple $F$ of nearby cycles complexes $\psi_F\sP$ for perverse sheaves $\sP$ on $X$ has been recently established in \cite{Wu}. This holds in particular for $\psi_F \bC_X$. To achieve the correspondence one has to
pay a price, for example for $\psi_F \bC_X$ one has to consider all $\bZ^r$ translates of the irreducible components of codimension 1 of the zero locus of $B_F$. Nevertheless analogs of Theorem \ref{thm3} $(i)$ and $(iii)$ have been shown for $\psi_F\sP$ for all perverse $\sP$ in \cite[Corollary 1.3, Theorem 1.4, Corollary 1.5]{Wu}.

\subsection{Monodromy conjecture for $F$.}
Conjecture \ref{conj1} $(i)$ with all $f_j$ elements of a multiplicatively - closed class of functions $\sC$ was shown in \cite{B} to be equivalent to its classical version for $r=1$ for functions in $\sC$. Thus Conjecture \ref{conj1} $(i)$ holds when $\sC$ is the class of plane curves and when $\sC$ is the class of hyperplane arrangements, by using results of Loeser and, respectively, Budur-Musta\c{t}\u{a}-Teitler.

One hopes to say the same for Conjecture \ref{conj1} $(ii)$, but only limited results are available in this direction. 

In the case of hyperplane arrangements, Conjecture \ref{conj1} $(ii)$ is known for $F=(f_1,\ldots,f_r)$ given by: any factorization $f=\prod_{j=1}^rf_j$ of a tame hyperplane arrangement $f$ by \cite{Ba}, generalizing a result for $r=1$ of \cite{W}; and for  the complete factorization of $f=\prod_{j=1}^rf_j$ of any arrangement by \cite{Wu1}, and, later by \cite[Example 4.9]{SV} using tropicalizations.

In \cite{BV} it is shown that Conjecture \ref{conj1} holds when one of the polynomials satisfies some genericity with respect to the other $f_j$. Fixing a log resolution of a product $f=f_1\ldots f_p$ of polynomials, one can define a class of log (very) generic polynomials $g$ using images of hypersurfaces defined by high powers of ample line bundles (respectively, of line bundles  whose classes sit in some subcone of the ample cone satisfying a genericity condition) on a compactification of the resolution. Then Conjecture \ref{conj1} $(i)$
 holds for the tuple $(f_1,\ldots,f_p,g)$ and for the product $fg$, if $g$ is  log generic. Conjecture \ref{conj1} $(ii)$ holds for the tuple $(f_1,\ldots,f_p,g)$ and for the product $fg$, if $g$ is log very-generic.

\medskip\noindent{\bf Acknowledgement.}  N.B. was supported by a Methusalem grant  and by G0B3123N from FWO. L.W. was supported by a start-up grant from Zhejiang University.

\end{document}